\documentclass{article}
\usepackage{maa-monthly}

\usepackage[capitalize,nameinlink]{cleveref}
\crefname{lemma}{lemma}{lemmas}
\crefname{equation}{}{equations}
\crefname{chapter}{Appendix}{chapters}
\crefname{item}{}{items}
\crefname{enumi}{}{}
\crefname{figure}{Figure}{Figures}

\newcommand{\scal}[2]{\left\langle{#1},{#2}  \right\rangle}
\newcommand{\nnn}{\ensuremath{{n\in{\mathbb N}}}}
\newcommand{\Id}{\ensuremath{\operatorname{Id}}}
\newcommand{\Fix}{\ensuremath{\operatorname{Fix}}}
\newcommand{\RR}{\ensuremath{\mathbb R}}
\newcommand{\NN}{\ensuremath{\mathbb N}}
\newcommand{\ZZ}{\ensuremath{\mathbb Z}}
\newcommand{\thalb}{\ensuremath{\tfrac{1}{2}}}

\theoremstyle{theorem}
\newtheorem{theorem}{Theorem}

\theoremstyle{definition}

\begin{document}

\title{Edelstein's Astonishing Affine Isometry}
\markright{Edelstein's Astonishing Affine Isometry}
\author{Heinz H.\ Bauschke, Sylvain Gretchko, \\Walaa M.\ Moursi, and Matthew Saurette}

\maketitle

\begin{abstract}
In 1964, Michael Edelstein presented an amazing affine isometry 
acting on the space of square-summable sequences. This operator has no fixed points,
but a suborbit that converges to $0$ while another escapes in norm to
infinity! We revisit, extend and sharpen his construction. 
Moreover, we sketch a connection to modern optimization and monotone operator theory.
\end{abstract}

\section{Introduction.}

Suppose that 
\begin{equation*}
\text{$X$ is a real Hilbert space,}
\end{equation*}
with inner product $\scal{\cdot}{\cdot}$ and induced norm
$\|\cdot\|$.
Assume that $T\colon X\to X$.
If there exists $0\leq \kappa<1$ and 
$\|Tx-Ty\|\leq\kappa\|x-y\|$, 
then the celebrated 
\emph{Banach contraction mapping principle} guarantees
that $T$ has a unique fixed point $\bar{x}=T\bar{x}$
and that no matter how the starting point $x_0\in X$ is chosen,
the sequence of iterates $(T^nx_0)_\nnn$ converges to $\bar{x}$.
But what is the situation when $T$ is merely \emph{nonexpansive}, i.e., 
\begin{equation*}
\|Tx-Ty\|\leq\|x-y\| 
\end{equation*}
for all $x,y$ in $X$?
Well, for starters, $T$ need not have a fixed point (consider translations).
Or, $T$ could have many fixed points (consider $\Id$, the identity).
And even if $T$ has fixed points, the iterates do not necessarily converge
to a fixed point (consider $-\Id$). 
We refer the reader to \cite{BC2017,GK,GR} for more on this subject. 
While these negative examples suggest complications, there is indeed more that
can be said and also a very interesting history!

It starts with Browder and Petryshyn who proved in 1966 \cite[Theorem~1]{BP} that 
\begin{equation*}
\Fix T \neq \varnothing
\;\;\Leftrightarrow\;\;
\text{for every/some $x_0\in X$, the sequence $(T^nx_0)_\nnn$ is bounded.}
\end{equation*}
(If one orbit is bounded, then so are all by the nonexpansiveness of $T$.)
Negating this, we obtain 
\begin{equation*}
\Fix T = \varnothing
\;\;\Leftrightarrow\;\;
\text{for every/some $x_0\in X$, the sequence $(\|T^nx_0\|)_\nnn$ is unbounded.}
\end{equation*}
The latter situation was carefully examined by Pazy in 1971 in \cite{Pazy71}.
That paper contains many fine results on the topic; however, it also states the 
stronger result that $\Fix T = \varnothing$ 
if and only if for every/some $x_0 \in X$, $\|T^nx_0\|\to\infty$; however, 
\emph{the proof is not convincing}. Indeed, in a later paper from 1977, 
Pazy states (see \cite[end of Section~1]{Pazy77}) that 
he does not know whether there exists a fixed-point-free nonexpansive map
with $\varliminf_{n\to\infty}\|T^nx_0\|<\infty$. 
In hindsight, a proof was impossible to obtain. 
Indeed, in 1964 --- even well before Pazy's first paper appeared ---
Edelstein in \cite{Edel64} 
constructed a fixed-point free affine isometry for which a subsequence 
converges even to $0$! 
On the positive side, Roehrig and Sine did prove in 1981 (see \cite[Theorem~2]{RS}) that an Edelstein-like example
is impossible in finite-dimensional Hilbert spaces, i.e., Euclidean spaces: 
\begin{multline*}
\Fix T = \varnothing
\;\;\Leftrightarrow\;\;
\text{$\|T^nx_0\|\to\infty$ for every/some $x_0\in X$}\\
\text{ \emph{provided $X$ is finite-dimensional.}}
\end{multline*}
The goal of this paper is to bring the amazing discrete dynamical system
discovered by Edelstein to a broader audience. It turns out that the material
is essentially accessible to undergraduate mathematics students.  
We generalize his result, 
provide full details, 
and present a new much smaller subsequence that also ``blows up to infinity.''
And last but not least, we interpret this example through the lens of splitting methods
and monotone operator theory! 
The paper is organized as follows.
In \cref{sec:affrot}, we review affine rotations which form the basic
building blocks of Edelstein's example.
These rotations are lifted to $\ell^2$ in \cref{sec:ell2}.
\cref{sec:estimating} deals with estimating the norm of the orbit starting
at $0$. A suborbit converging to $0$ is presented in \cref{sec:0suborbit}, 
while \cref{sec:inftysuborbit} provides a suborbit blowing up in norm to $\infty$!
In \cref{sec:DR}, we sketch a connection to the Douglas--Rachford splitting 
operator. We conclude with an epilogue in \cref{sec:epilogue}. 

\section{Affine rotations in $\RR^2$.}

\label{sec:affrot}
We start with the linear rotation matrix and a fixed vector, i.e., 
\begin{equation}
\label{e:defL}
L := L_\theta := 
\begin{pmatrix}
  \cos\theta & -\sin\theta\\
  \sin\theta & \cos\theta
\end{pmatrix} 
\text{~and~}
v = 
\begin{pmatrix}
  v_1\\
  v_2
\end{pmatrix}. 
\end{equation}
The affine rotation operator we will first investigate is 
\begin{equation}
\label{e:April30c}
Rx := R_\theta x := L_\theta x + v.
\end{equation}
Note that $R$ is an \emph{affine isometry} (also known as 
\emph{isometric affine mapping}), i.e., 
\begin{equation}
\label{e:April30a}
\|Rx-Ry\|=\|x-y\| \quad
\text{for all $x$ and $y$ in $\RR^2$.}
\end{equation}
(Also, $L$ is a linear isometry.)
Let's determine the fixed points of $R$:
for $x\in\RR^2$, we have 
\begin{equation}
\label{e:April29a}
\text{
$x\in\Fix R$
$\Leftrightarrow$
$x=Rx$
$\Leftrightarrow$
$x=Lx+v$
$\Leftrightarrow$
$(\Id-L)x=v$.
}
\end{equation}
Now 
\begin{equation*}
\Id-L = 
\begin{pmatrix}
  1-\cos\theta & \sin\theta\\
  -\sin\theta & 1-\cos\theta
\end{pmatrix}
\end{equation*}
has determinant
$(1-\cos\theta)^2 +\sin^2\theta= 2(1-\cos \theta )$.
Assuming that $\cos\theta<1$, the 
matrix $\Id-L$ is invertible with 
\begin{align*}
(\Id-L)^{-1} &= 
\frac{1}{2(1-\cos\theta)}\begin{pmatrix}
  1-\cos\theta & -\sin\theta\\
  \sin\theta & 1-\cos\theta
\end{pmatrix} \\
&= 
\frac{1}{2}\begin{pmatrix}
  1 & -\cot(\theta/2) \\
  \cot(\theta/2) & 1
\end{pmatrix};
\end{align*}
thus, also using \cref{e:April29a},
\begin{subequations}
\label{e:April29b}
\begin{align}
f &:= (\Id-L)^{-1}v 
= 
\frac{1}{2(1-\cos\theta)}
\begin{pmatrix}
(1-\cos\theta)v_1 - (\sin\theta)v_2\\
(\sin\theta)v_1+(1-\cos\theta)v_2
\end{pmatrix}\\
&= 
\frac{1}{2}\begin{pmatrix}
v_1-v_2\cot(\theta/2)\\
v_2+v_1\cot(\theta/2)
\end{pmatrix}
\end{align}
\end{subequations}
is the \emph{unique} fixed point of $R$ 
and 
\begin{equation*}
f-Lf=v.
\end{equation*}
(If $\cos\theta=1$, then $Rx=x+v$ and thus 
either (i) $v=0$ and $\Fix R = \RR^2$ or 
(ii) $v\neq 0$ and $\Fix R = \varnothing$.) 
It is now easily shown by induction that 
$L_\theta^n = L_{n\theta}$ and that 
\begin{equation}
\label{e:Rcompo}
R^nx = f+L_{n\theta}(x-f)
\quad 
\text{
for every $n\in\NN$.
}
\end{equation}
Let us specialize this further by setting 
\begin{equation}
\label{e:April29c}
v = \begin{pmatrix}
  v_1 \\
  v_2
\end{pmatrix} 
= 
\xi\cdot \begin{pmatrix}
  1-\cos\theta \\
  -\sin\theta
\end{pmatrix}, 
\end{equation}
where $\xi>0$ is a parameter.
(From now on, we can again include the case when $\cos\theta = 1$!)
With this particular assignment, the fixed point $f$ 
of $R$ found in \cref{e:April29b} simplifies to 
\begin{equation}
\label{e:R2f}
f = \begin{pmatrix}
  \xi \\
  0
\end{pmatrix}. 
\end{equation}
It follows from \cref{e:Rcompo} and \cref{e:defL} that 
\begin{equation*}
R^nx = \begin{pmatrix}
  \xi+(x_1-\xi)\cos(n\theta) -x_2\sin(n\theta)\\
  (x_1-\xi)\sin(n\theta) + x_2\cos(n\theta)
\end{pmatrix}; 
\end{equation*}
in particular, 
\begin{equation*}
R^n 0 = \xi\cdot \begin{pmatrix}
  1-\cos(n\theta) \\
  -\sin(n\theta) 
\end{pmatrix}. 
\end{equation*}
Using 
the half-angle formula for the squared sine, we obtain 
\begin{equation}
\label{e:normR0}
\|R^n0\|^2 = \xi^2\big(2-2\cos(n\theta)\big)
= 4 \xi^2\sin^2(n\theta/2). 
\end{equation}
More generally, tackling $R^nx$, we have
\begin{align*}
\|R^nx\|^2 &= \|x\|^2 + 2\xi x_1(\cos(n\theta)-1)-2\xi x_2\sin(n\theta) + 2\xi^2(1-\cos(n\theta))\\
  &= 
 \|x\|^2 - 4\xi x_1\sin^2(n\theta/2)\\
 &\quad -4\xi x_2\sin(n\theta/2)\cos(n\theta/2) + 4\xi^2\sin^2(n\theta/2). 
\end{align*}
Soon, we will ``lift'' $R$ from $\RR^2$ to $\ell^2$.
To do so, 
we develop some bounds. 
For $\alpha$ and $\beta$ in $\RR$, we clearly have 
$-(\alpha^2+\beta^2) \leq 2\alpha\beta\leq (\alpha^2+\beta^2)$; thus, 
$-2(\alpha^2+\beta^2) \leq -4\alpha\beta \leq 2(\alpha^2+\beta^2)$.
Hence
\begin{align*}
&\hspace{-1 cm} - 4\xi x_1\sin^2(n\theta/2)-4\xi x_2\sin(n\theta/2)\cos(n\theta/2) \\
&\leq 2(x_1^2 + \xi^2\sin^4(n\theta/2)) + 2(x_2^2+\xi^2\sin^2(n\theta/2)\cos^2(n\theta/2)) \\
&=  2\|x\|^2+2\xi^2\sin^2(n\theta/2)(\sin^2(n\theta/2)+\cos^2(n\theta/2)) \\
&= 2\|x\|^2+2\xi^2\sin^2(n\theta/2).
\end{align*} 
Similarly, 
\begin{equation*}
- 4\xi x_1\sin^2(n\theta/2)-4\xi x_2\sin(n\theta/2)\cos(n\theta/2) 
 \geq -2\|x\|^2-2\xi^2\sin^2(n\theta/2).
\end{equation*}
Altogether, we finally have 
\begin{equation}
\label{e:isometrybound}
-\|x\|^2 + 2\xi^2\sin^2(n\theta/2) \leq \|R^nx\|^2 \leq 3\|x\|^2 + 6\xi^2\sin^2(n\theta/2). 
\end{equation}

\section{From $\RR^2$ to $\ell^2$: lifting $R$ to $\mathbf{R}$.}

\label{sec:ell2}

From now on, $\ell^2$ is the real Hilbert space of all square summable sequences.
We think of $\ell^2$ here as the subset 
\begin{equation*}
\ell^2 \subsetneq \RR^2\times\RR^2\times\cdots 
\end{equation*}
and we will think of the $k$th plane as indexed by $k\in \{1,2,\ldots\}$.
Now consider a \emph{sequence} of angles $(\theta_k)_{k\geq 1}$, 
where $\theta_k$ is the angle for  $k$th plane, 
as well as a \emph{sequence} of positive parameters $(\xi_k)_{k\geq 1}$.
Set, in the spirit of \cref{e:April29c}, 
\begin{equation*}
v_k = \xi_k\cdot
\begin{pmatrix}
1-\cos(\theta_k)\\
-\sin(\theta_k)
\end{pmatrix} \in \RR^2
\end{equation*}
for $k\geq 1$.
Using \cref{e:normR0}, 
we estimate
\begin{equation}
\label{e:minbound}
\|v_k\|^2 = 4\xi_k^2\sin^2(\theta_k/2) \leq \min\{4\xi_k^2,\xi_k^2\theta_k^2\}. 
\end{equation}
Then $\mathbf{v} := (v_k)_{k\geq 1}$ will lie in $\ell^2$, e.g., 
if $\boldsymbol{\xi} := (\xi_k)_{k\geq 1}$ lies in $\ell^2$ or 
if $\boldsymbol{\xi}$ is bounded and $\boldsymbol{\theta} := 
(\theta_k)_{k\geq 1}\in\ell^2$. 
An ingenious choice by Edelstein \cite{Edel64} will lead to
a nice analysis: 
we set 
\begin{equation}
\label{e:Edelangles}
\theta_k = \frac{2\pi}{k!} \in (0,2\pi],
\end{equation}
which gives 
\begin{equation}
\label{e:Edelnormv}
\|\mathbf{v}\|^2 = \|(v_k)_{k\geq 1}\|^2 = 4\sum_{k\geq 1}\xi_k^2\sin^2(\pi/k!).
\end{equation}
We assume from now on that 
\begin{equation*}
\text{$(\xi_k)_{k\geq 1}$ is positive and decreasing }
\end{equation*}
but not necessarily strictly decreasing 
(indeed, Edelstein used $\xi_k\equiv 1$).
Using \cref{e:minbound}, we have 
\begin{equation}
\label{e:April29d}
\|\mathbf{v}\|^2 \leq 4\sum_{k \geq 1 }\xi_k^2\frac{\pi^2}{(k!)^2}
<4\pi^2\xi_1^2\sum_{k\geq 1}\frac{1}{k!} = 4\pi^2\xi_1^2(\exp(1)-1). 
\end{equation}
We are now ready to extend $R$ to the countable Cartesian product of 
Euclidean planes:
Set 
\begin{equation*}
\mathbf{R}\colon 
\mathbf{x}= (x_1,x_2,\ldots,x_k,\ldots)\mapsto (R_{\theta_1}x_1,R_{\theta_2}x_2,\ldots, 
R_{\theta_k}x_k,\ldots),
\end{equation*}
where each $x_k\in\RR^2$. 
From \cref{e:isometrybound}, \cref{e:Edelnormv} and \cref{e:April29d}, 
we have, for $\mathbf{x}\in\ell^2$, 
\begin{equation}
\label{e:April30e}
\|\mathbf{R}\mathbf{x}\|^2 \leq 
3 \|\mathbf{x}\|^2 + 6\sum_{k}\xi_k^2\sin^2(\pi/k!) 
= 3 \|\mathbf{x}\|^2 + (3/2)\|\mathbf{v}\|^2 < \infty
\end{equation}
and thus $\mathbf{R}$ is an affine operator from $\ell^2$ to itself!
Owing to \cref{e:April30a}, we have that 
\begin{equation*}
\|\mathbf{R}\mathbf{x}-\mathbf{R}\mathbf{y}\|
=\|\mathbf{x}-\mathbf{y}\| \quad
\text{for all $\mathbf{x}$ and $\mathbf{y}$ in $\ell^2$.}
\end{equation*}
What about the corresponding fixed point? 
Well, in view of \cref{e:R2f}, an (algebraic) fixed point of $\mathbf{R}$
is 
\begin{equation*}
\mathbf{f} := (\xi_1,0,\xi_2,0,\xi_3,0,\ldots)
\end{equation*} 
and 
\begin{equation*}
\|\mathbf{f}\|^2 = \sum_{k\geq 1}\xi_k^2 \in [0,+\infty].
\end{equation*}
So the original Edelstein choice 
$\xi_k\equiv 1$ leads to 
$\mathbf{f}\notin \ell_2$, i.e., 
\begin{equation*}
\Fix\mathbf{R} = \varnothing
\quad\text{provided that $\xi_k\equiv 1$.}
\end{equation*}
(Other choices with the same outcome are possible, e.g., 
consider $\xi_k \equiv 1/\sqrt{k}$.)

\section{Estimating $\|\mathbf{R}^n(0)\|$.}
\label{sec:estimating}

\begin{figure}
  \centering
   \includegraphics[width=0.7\linewidth]{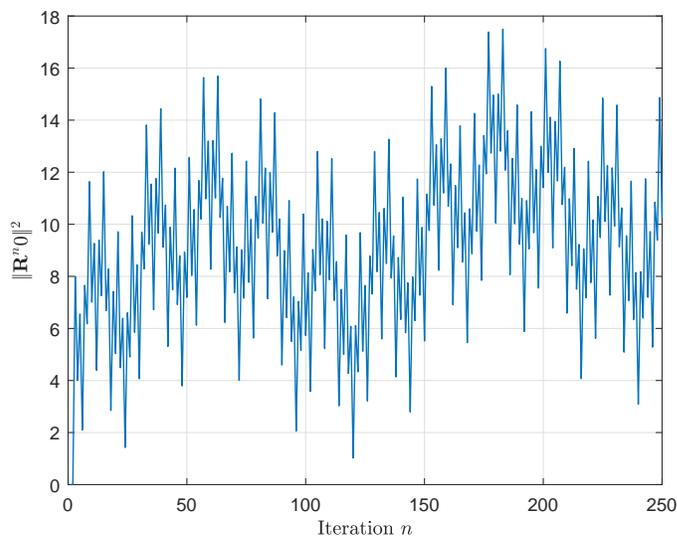}
   \caption{$\|\mathbf{R}^n0\|^2$ for $n\leq 250$.}
   \label{fig:250}
\end{figure}

\begin{figure}
  \centering
   \includegraphics[width=0.7\linewidth]{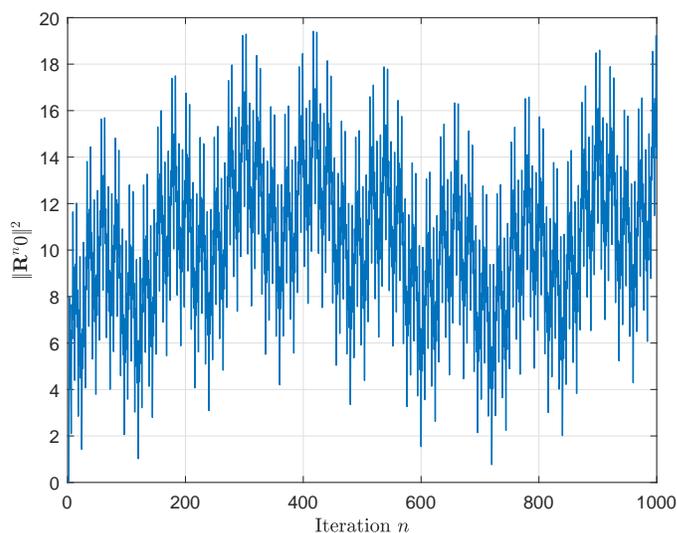}
   \caption{$\|\mathbf{R}^n0\|^2$ for $n\leq 1000$.}
   \label{fig:1000}
\end{figure}

Observe that 
\cref{e:isometrybound} yields
\begin{align*}
\|\mathbf{R}^n\mathbf{x}\|^2 
&\leq  3\|\mathbf{x}\|^2 + 6\sum_{k\geq 1}\xi_k^2\sin^2(n\pi/k!), 
\end{align*}
while \cref{e:normR0} gives the exact identity
\begin{equation}
\label{e:April27a}
\|\mathbf{R}^n0\|^2 = 4\sum_{k\geq 1}\xi_k^2\sin^2(n\pi/k!).
\end{equation}
For numerical computations, we need to estimate the convergence
of this infinite series. 
Indeed, 
\begin{align*}
0 &\leq \|\mathbf{R}^{n}0\|^2 - 
4\sum_{k=1}^{n}\xi_k^2\sin^2(\pi n/k!)
=
4\sum_{k\geq n+1}\xi_k^2\sin^2(\pi n/k!)\\
&\leq
4\pi^2\xi_1^2\sum_{k\geq n+1}(n/k!)^2\\
&\leq 
4\pi^2\xi_1^2\bigg(\frac{1}{(n+1)^2} + \frac{1}{(n+1)^2(n+2)^2}
+ \frac{1}{(n+1)^2(n+2)^2(n+3)^2}+\cdots\bigg)\\
&<
4\pi^2\xi_1^2\bigg(\frac{1}{(n+1)^2} + \frac{1}{(n+1)^4}
+ \frac{1}{(n+1)^6}+\cdots\bigg)\\
&=
\frac{4\pi^2\xi_1^2}{n(n+2)}\\
&\to 0 
\;\;\text{as $n\to \infty$.}
\end{align*}
Hence we will use for numerical computation 
\begin{equation*}
\|\mathbf{R}^{n}0\|^2 \approx 
4{\sum_{k=1}^{n}\xi_k^2\sin^2(\pi n/k!)}. 
\end{equation*}
We present $\|\mathbf{R}^n 0\|^2$ for the 
first 250 (respectively, 1000) iterates in 
\cref{fig:250} (respectively, \cref{fig:1000}). 
The source code is available at \cite{Sylvain}. 

Can you predict the long-term behavior from these plots?
We certainly would have guessed some form of periodic behavior.
However, this guess is far from the truth as we will see in the following sections. 
Indeed, a suborbit will converge to $0$ while another will blow up in norm to $\infty$. 
This bizarre behavior is not at all obvious from the plots!

\section{A suborbit that converges to $0$.}

\label{sec:0suborbit}

For any $n\geq 1$, we have 
\begin{align*}
&\|\mathbf{R}^{n!}0\|^2\\
&=
4\sum_{k\geq 1}\xi_k^2\sin^2(\pi n!/k!)
=
4\sum_{k\geq n+1}\xi_k^2\sin^2(\pi n!/k!)
\leq
4\pi^2\xi_1^2\sum_{k\geq n+1}(n!/k!)^2\\
&\leq
4\pi^2\xi_1^2\bigg(\frac{1}{(n+1)^2} + \frac{1}{(n+1)^2(n+2)^2}
+ \frac{1}{(n+1)^2(n+2)^2(n+3)^2}+\cdots\bigg)\\
&<
4\pi^2\xi_1^2\bigg(\frac{1}{(n+1)^2} + \frac{1}{(n+1)^4}
+ \frac{1}{(n+1)^6}+\cdots\bigg)\\
&=
\frac{4\pi^2\xi_1^2}{n(n+2)}\\
&\to 0 
\;\;\text{as $n\to \infty$.}
\end{align*}
Hence, we have the wonderfully weird suborbit result,
discovered first by Edelstein:
\begin{theorem} 
$\mathbf{R}^{n!}(0)\to 0$; in fact, 
\begin{equation*}
\|\mathbf{R}^{n!}(0)\| \leq O\big(\tfrac{1}{n}\big). 
\end{equation*}
\end{theorem}
This is very surprising and certainly not something 
easily deduced from the plots. 
(Of course, $n!$ grows really fast to $+\infty$.) 
However, in \cref{fig:250} and \cref{fig:1000}
we can indeed make out small values, namely for $\|\mathbf{R}^{120}0\|^2$ and 
for $\|\mathbf{R}^{720}0\|^2$.
This is not a coincidence since $120=5!$ and $720=6!$. 

In the next section, we will reveal a new suborbit that blows up
to infinity.

\section{A new suborbit that blows up to $\infty$.}

\label{sec:inftysuborbit}

Edelstein already provided a suborbit that goes to $\infty$ in norm.
Indeed, he suggested the sequence of indices
\begin{equation*} 
e_n := \tfrac{1}{2}\sum_{m=1}^n (n2^m)!
\end{equation*}
for which he proved that 
\begin{equation*}
\|\mathbf{R}^{e_n}0\|\to \infty.
\end{equation*}
It is hard to grasp just how fast $(e_n)_{n\geq 1}$ grows.
Indeed, the first three terms are 
\begin{equation*}
(e_1,e_2,e_3) = (1,20172,310224200866619959181160).
\end{equation*}
We won't list 
$e_4,e_5,e_6$ but point out that they have 
$89,285,828$ decimal digits, respectively. 
(The source code is available at \cite{Sylvain}.)
Instead, we present a 
new special subsequence that will achieve the same result
but with much smaller indices. 
Set 
\begin{equation*}
s_n := 1 + \sum_{m=1}^{n-1}\lceil m/2\rceil (m+2)!
\end{equation*}
for every $n\geq 1$. 
Of course, $s_n$ approaches infinity 
but much slower than $(e_n)_{n\geq 1}$;
indeed, 
\begin{equation*}
(s_1,s_2,s_3,s_4,s_5,s_6)
= (7,31,271,1711,16831,137791).
\end{equation*}
We now start to analyze the term $s_n$.
First, 
\begin{align*}
\frac{s_n}{k!} 
&= \frac{1}{k!} + \sum_{m=1}^{n-1}\lceil m/2\rceil (m+2)!/k!\\
&= \frac{1}{k!} + 
\sum_{m=1}^{k-3}\lceil m/2\rceil (m+2)!/k!
+ \sum_{m=k-2}^{n-1}\lceil m/2\rceil \underbrace{(m+2)!/k!}_{\in \ZZ},
\end{align*}
where $k-3 \leq n-1$, i.e., 
$k\leq n+2$.
Next, we will work on deriving bounds on the fractional part of 
$s_n/k!$ which will be used later. 
Assuming that $k-3\geq 1$, i.e., $k\geq 4$, we have 
\begin{align*}
\frac{1}{k!} + 
\sum_{m=1}^{k-3}\lceil m/2\rceil (m+2)!/k!
&> 
0+
\sum_{m=k-3}^{k-3}\lceil m/2\rceil (m+2)!/k!\\
&=\lceil (k-3)/2\rceil / k
\geq (k-3)/(2k)\\
&=\frac{1}{2}-\frac{3}{2k}.
\end{align*}
We now turn to an upper bound. 
Using $\lceil q/2 \rceil < (q+2)/2$
and requiring that $k\geq 5$, 
we obtain 
\begin{align*}
&\sum_{m=1}^{k-3}\frac{\lceil m/2\rceil (m+2)!}{k!}\\
&=
\frac{\lceil (k-3)/2\rceil (k-1)!}{k!} 
+ 
\frac{\lceil (k-4)/2\rceil (k-2)!}{k!} 
+ 
\sum_{m=1}^{k-5}\frac{\lceil m/2\rceil (m+2)!}{k!}\\
&=
\frac{\lceil (k-3)/2\rceil}{k} 
+ 
\frac{\lceil (k-4)/2\rceil }{k(k-1)} 
+ 
\sum_{m=1}^{k-5}\frac{\lceil m/2\rceil}{(m+3)\cdots(k-2)(k-1)k}\\
&< 
\frac{k-1}{2k} 
+ 
\frac{k-2}{2k(k-1)} 
+ 
\sum_{m=1}^{k-5}\frac{(m+2)/2}{(k-2)(k-1)k}\\
&= \frac{(k+1)(2k^2-7k+4)}{4(k-2)(k-1)k}.
\end{align*}
It follows that 
\begin{align*}
\frac{1}{k!} + \sum_{m=1}^{k-3}\frac{\lceil m/2\rceil (m+2)!}{k!}
&< 
\frac{1}{k!}
+ \frac{(k+1)(2k^2-7k+4)}{4(k-2)(k-1)k}\\
&< 
\frac{1}{k!}
+ \frac{(k+1)(2k^2-6k+4)}{4(k-2)(k-1)k}\\
&= 
\frac{1}{k!}
+ \frac{k+1}{2k}
= \frac{1}{2}+\frac{1}{k!}+\frac{1}{2k}\\
&<\frac{1}{2} + \frac{1}{(k)(4)(3)(2)} + \frac{1}{2k}\\
&= \frac{1}{2}+\frac{13}{24k}. 
\end{align*}
Altogether, the \emph{fractional part} $\{s_n/k!\}$ of 
$s_n/k!$ satisfies, for $k\geq 10$, 
\begin{equation*}
\frac{1}{3}<\frac{7}{20}=\frac{1}{2} - \frac{3}{20}\leq \frac{1}{2}-\frac{3}{2k}<\left\{\frac{s_n}{k!}\right\}< \frac{1}{2}+\frac{13}{24k}
\leq \frac{1}{2} + \frac{13}{240} = \frac{133}{240}<\frac{2}{3}. 
\end{equation*}
Combining this with the fact that $\sin^2$ over $[\pi/3,2\pi/3]$ has minimum 
value $3/4$ results in 
\begin{equation*}
\frac{3}{4}<\sin^2\Big(\frac{s_n}{k!}\pi\Big)
\quad 
\text{whenever $10\leq k \leq n+2$.}
\end{equation*}
Using \cref{e:April27a}, we thus get 
\begin{equation}
\label{e:April28a}
\|\mathbf{R}^{s_n}0\|^2 
= 
4\sum_{k=1}^\infty \xi_k^2\sin^2(s_n\pi/k!)
> 3\sum_{k=10}^{n+2} \xi_k^2
\geq 3(n-7)\xi_{n+2}^2.
\end{equation}

We now turn toward an upper bound.
Again using $\lceil q/2 \rceil < (q+2)/2$ 
as well as  an easy induction in \cref{e:easyind}, 
we obtain 
\begin{subequations}
\begin{align}
s_n 
&=
1 + \sum_{m=1}^{n-1}\lceil m/2\rceil (m+2)!
< 
1+ \sum_{m=1}^{n-1}\thalb(m+2)(m+2)!\\
&=
1+\thalb(n+2)!-3 \label{e:easyind}\\
&<
\thalb (n+2)!
\end{align}
\end{subequations}
which implies
\begin{align*}
\frac{s_n}{(n+3)!} &< \frac{1}{2(n+3)}, \;\;
\frac{s_n}{(n+4)!} < \frac{1}{2(n+3)(n+4)},\\
\frac{s_n}{(n+5)!} &< \frac{1}{2(n+3)(n+4)(n+5)}, \;\; \ldots \; . 
\end{align*}
Thus 
\begin{align*}
4\sum_{k=n+3}^\infty \xi_k^2\sin^2(s_n\pi/k!)
&<
\pi^2\sum_{k=n+3}^\infty \frac{\xi_k^2}{(n+3)^2(n+4)^2\cdots k^2}\\
&<
\frac{\pi^2\xi_{n+3}^2}{(n+3)^2}\bigg(1 + \frac{1}{(n+3)^2} + \frac{1}{(n+3)^4} + \cdots \bigg)\\
&= \frac{\pi^2\xi_{n+3}^2}{(n+3)^2-1}.
\end{align*}
Again using \cref{e:April27a}, we thus have 
\begin{equation}
\label{e:April28b}
\|\mathbf{R}^{s_n}0\|^2 = 
4\sum_{k=1}^\infty \xi_k^2\sin^2(s_n\pi/k!)
<
4\sum_{k=1}^{n+2}\xi_k^2 + 
\frac{\pi^2\xi_{n+3}^2}{(n+3)^2-1}.
\end{equation}

Combining \cref{e:April28a} with \cref{e:April28b}, we obtain 
the following result, which summarizes our work in this section:

\begin{theorem}
If $n\geq 8$, then 
\begin{align*}
3(n-7)\xi_{n+2}^2
&\leq 3\sum_{k=10}^{n+2}\xi_k^2\\
&\leq \|\mathbf{R}^{s_n}0\|^2
<  
4\sum_{k=1}^{n+2}\xi_k^2 + 
\frac{\pi^2\xi_{n+3}^2}{(n+3)^2-1}\\
&\leq 4(n+2)\xi_1^2 + \frac{\pi^2\xi_{n+3}^2}{(n+3)^2-1}.
\end{align*}
Consequently, we have the implications
\begin{equation*}
\lim_{k\to\infty}\xi_k > 0
\;\;\Rightarrow\;\;
(\xi_k)_{k\geq 1}\notin\ell_2
\;\;\Rightarrow\;\;
\|\mathbf{R}^{s_n}0\|\to\infty.
\end{equation*}
and 
\begin{equation*}
\lim_{k\to\infty}\xi_k > 0
\;\;\Rightarrow\;\;
\|\mathbf{R}^{s_n}0\| = O\big(\sqrt{n}\,\big). 
\end{equation*}
\end{theorem}

Note that this covers the Edelstein set-up where $\xi_k\equiv 1$, 
albeit with the much smaller sequence of indices 
$(s_{n})_{n\geq 1}$ compared to $(e_n)_{n \geq 1}$.

\section{Optimization and the Douglas--Rachford algorithm.}

\label{sec:DR}

In this section, we place the example into a different framework.
Indeed, we have so far studied the world of nonexpansive mappings.
There are two related worlds, namely the ones featuring 
firmly nonexpansive mappings and monotone operators.
We start with \emph{firmly nonexpansive} mappings, which are mappings
on $X$ that satisfy
\begin{equation*}
\|Tx-Ty\|^2 \leq \scal{Tx-Ty}{x-y} 
\quad\text{for all $x$ and $y$ in $X$.}
\end{equation*}
It is not hard to show that $T$ is firmly nonexpansive if
and only if $2T-\Id$ is nonexpansive, which implies that 
the firmly nonexpansive mappings are precisely those that
can be written as the average of the identity and nonexpansive mappings
(see, e.g., \cite{BC2017} for more on this).

Let's identify the firmly nonexpansive counterpart of 
the operator $R=R_\theta$ considered in \cref{e:April30c}, 
with $0<\theta\leq 2\pi$. 
We thus set 
\begin{equation*}
T := T_\theta := \tfrac{1}{2}\Id + \tfrac{1}{2}R_\theta.
\end{equation*}
Let $x\in X$. Using the notation and results of \cref{sec:affrot}, it can be 
shown by induction that 
\begin{equation*}
T^nx = f + \frac{1}{2^n}(\Id+L)^n(x-f). 
\end{equation*}
Using the double-angle formula for sine and half-angle formula for cosine 
(with angle $\theta/2$), we obtain
$\Id+L_\theta = 2\cos(\theta/2)L_{\theta/2}$.
Hence $(\Id+L_{\theta})/2 = \cos(\theta/2)L_{\theta/2}$
and thus $\big((\Id+L_{\theta})/2 \big)^n = \cos^n(\theta/2)L_{n\theta/2}$.
It follows that if $\theta<2\pi$, then 
\begin{equation}
\label{e:Tandcos}
T^nx = f + \cos^n(\theta/2)L_{n\theta/2}(x-f) \to f
\end{equation}
as predicted by  Rockafellar's proximal point algorithm 
(see, e.g., \cite{Rockprox}); 
indeed, 
\begin{equation*}
\|T^nx-f\| = \cos^n(\theta/2)\|x-f\| \to 0
\end{equation*}
and $T^nx\to f$
with a sharp linear rate of $\cos(\theta/2)$.
And if $\theta=2\pi$, then $T_{2\pi}=\Id$ and we have immediate
and even finite convergence to a fixed point!

The operator $T$ on $\RR^2$ gives rise to an induced operator 
$\mathbf{T}\colon \ell^2\to\ell^2$ (because $\mathbf{R}$ is 
well-defined on $\ell^2$ by \eqref{e:April30e} and hence so 
is $\mathbf{T}=(\Id+\mathbf{R})/2$). 
While $T$ does have the unique fixed point $f$ (unless $\theta=2\pi$)
the same is no longer true in general
for $\mathbf{T}$ because an algebraic fixed point 
$\mathbf{f}$ may fail to lie in $\ell^2$.
($\Fix \mathbf{T}$ cannot be a singleton because 
$\Fix T_{\theta_1} = \Fix T_{2\pi} = \RR^2$.)

Let us now interpret this from an optimization/feasibility perspective.
Let $U$ and $V$ be nonempty closed convex subsets of $\RR^2$ with
$U\cap V\neq\varnothing$, and denote the 
projection mappings onto $U,V$ by $P_U,P_V$, respectively.
Then the associated Douglas--Rachford splitting operator is 
\begin{equation*}
D := D_{V,U} = \Id-P_U+P_V(2P_U-\Id).
\end{equation*}
It is known that the sequence $(P_UD^nx)_{n\geq 1}$
converges to some point in $U\cap V$; this is a well-known
method to solve feasibility (and even optimization) problems.
Now suppose that 
\begin{equation*}
U = \RR(1,0)=\RR\times\{0\}
\;\;\text{and}\;\;
V = f + \RR(\cos(\theta/2),\sin(\theta/2)),
\end{equation*}
where $f=(\xi,0)$ is as in \cref{e:R2f}.
Clearly, $U\cap V =\{f\}$ (unless $\theta=2\pi$, in which case
$U\cap V=U=V$).
Using linear algebra, one may check that 
\begin{equation*}
D = T;
\end{equation*}
put differently, the Douglas--Rachford operator for the feasibility problem
of finding a point in $U\cap V$ is exactly the firmly nonexpansive counterpart
of the Edelstein affine isometry!
Working towards the $\ell_2$ version, 
consider first the pure Cartesian products
\begin{equation*}
\widetilde{\mathbf{U}} = (\RR\times\{0\})\times(\RR\times\{0\})\times \cdots 
\end{equation*}
and 
\begin{align*}
\widetilde{\mathbf{V}} 
&= \big((\xi_1,0)+\RR(\cos(\theta_1/2),\sin(\theta_1/2))\big) 
\times \big(\xi_2,0)+\RR(\cos(\theta_2/2),\sin(\theta_2/2))\big) \\
&\quad \times  \cdots\\
&= \big((\xi_1,0)+\RR(\cos(\pi/1!),\sin(\pi/1!))\big) 
\times \big(\xi_2,0)+\RR(\cos(\pi/2!),\sin(\pi/2!))\big) \\
&\quad \times  \cdots\\
&= \big(\RR\times\{0\}\big)
\times \big(\{\xi_2\}\times\RR)\big) 
\times \cdots \times \big((\xi_k,0)+\RR(\cos(\pi/k!),\sin(\pi/k!))\big) \\
&\quad \times \cdots 
\end{align*}
for which 
\begin{equation*}
\widetilde{\mathbf{U}}\cap \widetilde{\mathbf{V}}
= \big(\RR\times\{0\})\times
\{(\xi_2,0)\} 
\times \cdots \times
\{(\xi_k,0)\} \times\cdots \;. 
\end{equation*}
Now set 
\begin{equation*}
\mathbf{U} = \widetilde{\mathbf{U}} \cap \ell^2
\;\;\text{and}\;\;
\mathbf{V} = \widetilde{\mathbf{V}} \cap \ell^2. 
\end{equation*}
Then our $\mathbf{T}$ is precisely the Douglas--Rachford operator 
for finding a point in $\mathbf{U}\cap \mathbf{V}$.
Note that this set is possibly empty when $\boldsymbol{\xi}\notin\ell^2$, 
as is the case for Edelstein's original $\xi_k\equiv 1$. 

We finally turn to the third world, the world of (maximally) monotone operators.
The unique maximally monotone operator associated with $R$ and with $T$ is 
\begin{equation*}
M := M_\theta : =T_\theta^{-1}-\Id.
\end{equation*}
Let us find out what $M$ is. 
We start by 
inverting $T$, for which \cref{e:Tandcos} surely helps.
So we write $x=Ty$. 
Let us abbreviate $c := \cos(\theta/2)$
and $K := L_{\theta/2}$. 
Assume that $c\neq 0$, i.e., $\theta\neq\pi$. 
Note that $K^{-1} = L_{-\theta/2} = L_{\theta/2}^* = K^*$. 
Then
$x=Ty=f+c K(y-f)$
if and only if
$c^{-1}(x-f) = K(y-f)$
if and only if
$c^{-1}K^{-1}(x-f) = y -f$
if and only if
$f+c^{-1}K^*(x-f) = y$. 
Hence $Mx = f+c^{-1}K^*(x-f)-x = (c^{-1}K^*-\Id)(x-f)$. 
Switching back to the original notation, we have
\begin{align*}
Mx &= \left(\frac{1}{\cos(\theta/2)} 
\begin{pmatrix}
\cos(-\theta/2) & -\sin(-\theta/2)\\
\sin(-\theta/2) & \cos(-\theta/2)
\end{pmatrix}
- 
\begin{pmatrix}
1 & 0 \\
0 & 1
\end{pmatrix}
\right)(x-f)\\
&= 
\begin{pmatrix}
0 & \tan(\theta/2)\\
-\tan(\theta/2) & 0 
\end{pmatrix}
(x-f)\\
&= \tan(\theta/2)
\begin{pmatrix}
0 & 1\\
-1 & 0 
\end{pmatrix}
(x-f).
\end{align*}
If $\theta = \pi$, then $Tx\equiv f$ and hence 
the set-valued inverse of $T$ maps $f$ to $\RR^2$ and anything
else to the empty set.
To summarize,
\begin{equation}
\label{e:M}
M_\theta x = 
\begin{cases}
\tan(\theta/2)
\begin{pmatrix}
0 & 1\\
-1 & 0 
\end{pmatrix}(x-f), &\text{if $\theta\neq\pi$;}\\
\RR^2, &\text{if $\theta=\pi$ and $x=f$;}\\
\varnothing, &\text{if $\theta=\pi$ and $x\neq f$.}
\end{cases}
\end{equation}
Monotone operator theory predicts that $M$ is \emph{monotone}, i.e., 
$\scal{x^*-y^*}{x-y} \geq 0$ 
for all $x^*\in Mx$ and $y^*\in My$. 
Indeed, our operator $M$ satisfies this inequality even as equality.

Recall that the Edelstein angles given by \cref{e:Edelangles} are
$\theta_1 = 2\pi$, $\theta_2 = \pi$, $\theta_3 = \pi/3$, \ldots, 
$\theta_k=2\pi/k!$, \ldots\ . 
Let $k\geq 3$. Then $\theta_k\neq \pi$.
Writing $x=(x_1,x_2)$, $t = \tan(\theta_k/2)$, and 
$\xi=\xi_k$, it follows from \cref{e:M} and \cref{e:R2f} that 
$Mx = t(-x_2,x_1-\xi)$ and thus 
\begin{align*}
\|Mx\|^2 &= t^2(x_2^2+(x_1-\xi)^2)
= t^2(\|x\|^2 + \xi^2-2\xi x_1)\\
&\leq t^2(\|x\|^2 + 2\xi^2+x_1^2)
\leq 2t^2(\|x\|^2+\xi^2). 
\end{align*}
Switching to the product space version $\mathbf{M}$
acting on $\RR^2\times\RR^2\times \cdots$, 
we deduce that 
if ${\bf u} = (u_k)_{k\geq 1}\in \mathbf{M}\mathbf{x}$ and 
$\mathbf{x}\in \ell^2$, then 
\begin{align*}
\|\mathbf{u}\|^2
-\|u_1\|^2-\|u_2\|^2 
&= \sum_{k\geq 3}\|u_k\|^2
\leq 2\sum_{k\geq 3} \tan^2(\theta_k/2)(\|x_k\|^2+\xi_k^2)\\
&\leq 2(\|\mathbf{x}\|^2+\xi_1^2)\sum_{k\geq 3}\tan^2(\pi/k!)<\infty
\end{align*}
because 
\begin{align*}
\tan^2(\pi/k!)&= (\sin(\pi/k!)/\cos(\pi/k!))^2
\leq (\pi/k!)^2/\cos^2(\pi/3!)\\
&= \pi^2/(\sqrt{3}/2)^2/(k!)^2
=(4\pi^2/3)/(k!)^2
\end{align*}
and the comparison test applies. %
Hence for $\mathbf{x}=(x_1,x_2,\ldots)\in\ell^2\subsetneqq \RR^2\times\RR^2\times\ldots$ we have
\begin{equation*}
\mathbf{M}\mathbf{x}\neq\varnothing 
\quad\Leftrightarrow\quad
x_2=(\xi_2,0).
\end{equation*}

This concludes our journey featuring the Edelstein operator --- 
we hope you enjoyed the ride as much as we did!

\section{Epilogue.}

\label{sec:epilogue}

The authors believe that Jonathan (Jon) 
Borwein would have liked the material in this paper for several reasons:
While at Dalhousie University, Jon and Michael Edelstein were actually colleagues.
Jon enjoyed nonexpansive mappings and published extensively in this area (see,
e.g., \cite{BS84}). 
In his optimization work, Jon spent significant 
time of his last years working 
on the Douglas--Rachford algorithm (see, e.g., \cite{ABT13}). 
Last but not least, 
some of the numbers in this paper were obtained by computation 
(see \cite{Sylvain}) --- as a co-founder of experimental mathematics (see, 
e.g., \cite{BBCGLM}), Jon would have 
enjoyed the flavor of these results and the concreteness of the examples.

\begin{acknowledgment}{Acknowledgment.}
The authors thank the referees and the editors for 
their careful reading and constructive comments. 
HHB was partially supported by the Natural Sciences and
Engineering Research Council of Canada.
WMM was partially supported by a Natural Sciences and
Engineering Research Council of Canada Postdoctoral Fellowship.
\end{acknowledgment}

\begin{biog}
\item[Heinz Bauschke] 
is a former doctoral student of Jonathan Borwein and currently a Professor of
Mathematics at the University of British Columbia (Okanagan campus in
Kelowna, B.C., Canada). His main areas of interest are in convex analysis,
optimization, and monotone operator theory. He has published more than 100
papers and the book \emph{Convex Analysis and Monotone Operator Theory in Hilbert
Spaces}. 
\begin{affil}
Department of Mathematics, UBC Okanagan, Kelowna, B.C. V1V 1V7,
Canada\\
heinz.bauschke@ubc.ca
\end{affil}

\item[Sylvain Gretchko] 
is a Software Engineer and a graduate student in mathematics at the 
University of G\"ottingen. 
He is particularly interested in experimental mathematics.
\begin{affil}
Faculty of Mathematics and Computer Science,
University of G\"ottingen,
37073 G\"ottingen,
Germany\\
sylvain.gretchko@gmail.com
\end{affil}

\item[Walaa Moursi] 
is an Assistant Professor in the Department of Combinatorics and Optimization 
at the University of Waterloo. 
Her research interests are convex analysis, 
monotone operator theory and continuous optimization.
\begin{affil}
Department of Combinatorics and Optimization,
Faculty of Mathematics,
University of Waterloo,
Waterloo, Ontario N2L 3G1,
Canada\\
walaa.moursi@uwaterloo.ca
\end{affil}

\item[Matthew Saurette]
is an undergraduate student in mathematics at the University of British Columbia
(Okanagan campus in Kelowna, B.C., Canada). His interests are in optimization
and mathematical biology.
\begin{affil}
Department of Mathematics, UBC Okanagan, Kelowna,
B.C.~V1V 1V7, Canada\\ 
matthews14@hotmail.ca
\end{affil}

\end{biog}

\vfill\eject

\end{document}